\documentclass[preprint,12pt]{elsarticle}
\usepackage{amsmath}
\usepackage{amssymb}

\begin{document}
\begin{frontmatter}

\ead{aaalikhanov@gmail.com}

\title{ A new difference scheme for the time fractional diffusion equation}

\author{Anatoly A. Alikhanov}

\address{Kabardino-Balkarian State University, ul. Chernyshevskogo 173,
 Nalchik,  360004,   Russia}

\begin{abstract}{ In this paper we construct a new difference analog of the
Caputo fractional derivative  (called the $L2$-$1_\sigma$ formula).
The basic properties of this difference operator are investigated
and on its basis some difference schemes generating approximations
of the second and fourth order in space and the second order in time
for the time fractional diffusion equation with variable
coefficients are considered.  Stability of the suggested schemes and
also their convergence in the grid  $L_2$ - norm with the rate equal
to the order of the approximation error are proved. The obtained
results are supported by the numerical calculations carried out for
some test problems.}
\end{abstract}


\end{frontmatter}

%%%%%%%%%%%%%%%%%%%%%%%%%%%%%%%%%%%%%%%%%%%%%%%%%%%%%%%%%%%%%%%%%%%%%%%
\section{Introduction}
Recently a noticeable growth of the attention of researches to the
fractional differential equations has been observed. It is caused by
numerous effective applications of fractional calculation to various
areas of science and engineering
\cite{Nakh:03,Old,Podlub:99,Hilfer:00,Kilbas:06,Uchaikin:08}. For
example, mathematical language of fractional derivatives is
irreplaceable for the description of the physical process of
statistical transfer and, as it is known, leads to diffusion
equations of fractional orders \cite{Nigma,Chuk}.

Consider the time fractional diffusion equation with variable
coefficients
\begin{equation}\label{ur01}
\partial_{0t}^{\alpha}u(x,t)=\mathcal{L}u(x,t)+f(x,t) ,\,\, 0<x<l,\,\,
0<t\leq T,
\end{equation}
\begin{equation}
u(0,t)=0,\quad u(l,t)=0,\quad 0\leq t\leq T, \quad
u(x,0)=u_0(x),\quad 0\leq x\leq l,\label{ur02}
\end{equation}
where
\begin{equation}
\partial_{0t}^{\alpha}u(x,t)=\frac{1}{\Gamma(1-\alpha)}\int\limits_{0}^{t}\frac{\partial
u(x,\eta)}{\partial\eta}(t-\eta)^{-\alpha}d\eta, \quad 0<\alpha<1
\label{ur02.01}
\end{equation}
is the Caputo derivative of the order $\alpha$,
$$
\mathcal{L}u(x,t)=\frac{\partial }{\partial
x}\left(k(x,t)\frac{\partial u}{\partial x}\right)-q(x,t)u,
$$
$k(x,t)\geq c_1>0$, $q(x,t)\geq 0$ and $f(x,t)$ are sufficiently
smooth functions.

The time fractional diffusion equation represents a linear integro -
differential equation. Its solution not always can be found
analytically; therefore it is necessary to use numerical methods.
However, unlike the classical case, we require information about all
the previous time layers, when numerically approximating a time
fractional diffusion equation on a certain time layer. For that
reason algorithms for solving the time fractional diffusion
equations are rather time-consuming even in one - dimensional case.
Upon transition to two - dimensional and three - dimensional
problems their complexity considerably increases. In this regard
constructing stable differential schemes of higher order
approximation is a very important task.

A widespread difference approximation of fractional
derivative  (\ref{ur02.01}) is the so-called $L1$ method
\cite{Old,Sun} which is defined as follows
\begin{equation}
\partial_{0t_{j+1}}^{\alpha}u(x,t)=\frac{1}{\Gamma(1-\alpha)}\sum\limits_{s=0}^{j}\frac{u(x,t_{s+1})-u(x,t_{s})}{t_{s+1}-t_{s}}\int\limits_{t_{s}}^{t_{s+1}}
\frac{d\eta}{(t_{j+1}-\eta)^{\alpha}}+r^{j+1}, \label{ur02.02}
\end{equation}
where $0=t_0<t_1<\ldots<t_{j+1}$, and $r^{j+1}$ is the local truncation
error. In the case of the uniform mesh, $\tau=t_{s+1}-t_s$, for all
$s=0,1,\ldots, j+1$, it was proved that
$r^{j+1}=\mathcal{O}(\tau^{2-\alpha})$ \cite{Sun2,Lin,Alikh_arxiv3}.
The $L1$ method has been widely used for solving the fractional
differential equations with Caputo derivatives
\cite{Sun2,Lin,Alikh_arxiv3,ShkhTau:06,Liu:10,Alikh:12,Delic}.

Difference schemes of the increased order of approximation such as
the compact difference scheme \cite{Liu:10,Sun3,Sun4,Sun5} and
spectral method \cite{Lin,Lin2,Xu} were applied to improve the
spatial accuracy of fractional diffusion equations. However, it is
rather difficult to get a high-order time approximation due to the
singularity of fractional derivatives.

 A good approximation of the  $L1$ method is observed in case of a nonuniform mesh,
  when it is refined  in a neighborhood of the  point  $t_{j+1}$ \cite{Sun}.
Though the nonuniform mesh turns out to be more effective in
comparison with the uniform one, it will not generate the second
order of approximation in all points of the mesh.

In \cite{Sun6} a new difference analog of the Caputo fractional
derivative with the order of approximation
$\mathcal{O}(\tau^{3-\alpha})$, called $L1-2$ formula, is
constructed. On the basis of this formula calculations of difference
schemes for the time-fractional sub-diffusion equations in bounded
and unbounded spatial domains and the fractional ODEs are carried
out. If the stability and convergence of difference schemes from
\cite{Sun6} will be strictly proved, then this will undoubtedly be a
significant progress in computing the time-fractional partial
differential equations.

Using the energy inequality method, a priori estimates for the
solution of the Dirichlet and  Robin boundary value problems for the
diffusion-wave equation with Caputo fractional derivative have been
obtained in \cite{Alikh:12,Alikh:10}.

In this paper a new difference analog of the fractional Caputo
derivative with the order of approximation
$\mathcal{O}(\tau^{3-\alpha})$ for each  $\alpha\in(0,1)$  is
constructed. Properties of the obtained difference operator are
studied. Difference schemes of the second and fourth order of
approximation in space and the second order in time for the time
fractional diffusion equation with variable coefficients are
constructed. Using the method of energy inequalities, the stability
and convergence  of these schemes in the mesh $L_2$ - norm are
proved. Numerical calculations of some test problems confirming
reliability of the obtained results are carried out.

\section{Family of difference schemes. Stability and convergence}

In this section, families of difference schemes in a general form set
on a non-uniform time mesh are investigated. A criterion of the
stability of the difference schemes in the mesh $L_2$ - norm is
obtained. The convergence of solutions of the difference schemes to the
solution of the corresponding differential problem with the rate
equal to the order of the approximation error is proved.

\subsection{Family of difference schemes}

In the rectangle  $\overline Q_T=\{(x,t): 0\leq x\leq l,\, 0\leq
t\leq T\}$ we introduce the mesh $\overline
\omega_{h\tau}=\overline\omega_{h}\times\overline\omega_{\tau}$,
where $\overline\omega_{h}=\{x_i=ih, \, i=0, 1, \ldots, N;\,
hN=l\}$, $\overline\omega_{\tau}=\{t_j: \, 0=t_0<t_1<t_2<\ldots
<t_{M-1}<t_{M}=T\}$.

Basically the family of difference schemes, approximating problem
(\ref{ur01})--(\ref{ur02}) on the mesh $\overline \omega_{h\tau}$,
has the form

\begin{equation}\label{ur03}
{_g}\Delta_{0t_{j+1}}^{\alpha}y_i=\Lambda y^{(\sigma_{j+1})}_i
+\varphi_i^{j+1}, \quad i=1,2,\ldots,N-1,\quad j=0,1,\ldots,M-1,
\end{equation}
\begin{equation}
y(0,t)=0,\quad y(l,t)=0,\quad t\in \overline \omega_{\tau}, \quad
y(x,0)=u_0(x),\quad  x\in \overline \omega_{h},\label{ur03.1}
\end{equation}
where
\begin{equation}
{_g}\Delta_{0t_{j+1}}^{\alpha}y_i=\sum\limits_{s=0}^{j}\left(y_i^{s+1}-y_i^s\right)g_{s}^{j+1},\quad
g_{s}^{j+1}>0, \label{ur03.2}
\end{equation}
is a difference analog of the Caputo derivative of the order
$\alpha$ ($0<\alpha<1$), $\Lambda$ is  a difference operator
approximating the continuous operator  $\mathcal{L}$, such that the
operator $-\Lambda$ preserves its positive definiteness ($(-\Lambda
y,y)\geq \varkappa\|y\|^2$, $\varkappa>0$), for example
\begin{equation}\label{ur03.5}
(\Lambda y)_i=\left((ay_{\bar
x})_x-dy\right)_i=\frac{a_{i+1}y_{i+1}-(a_{i+1}+a_i)y_i+a_iy_{i-1}}{h^2}-d_iy_i,
\end{equation}
$a$, $d$ and $\varphi$ are the mesh functions approximating  $k$, $q$
and $f$, respectively,
$y^{(\sigma_{j+1})}=\sigma_{j+1}y^{j+1}+(1-\sigma_{j+1})y^{j}$,
$0\leq\sigma_{j+1}\leq1$, at $j=0,1,\ldots,M-1$, $y_{\bar
x,i}=(y_i-y_{i-1})/h$, $y_{x,i}=(y_{i+1}-y_{i})/h$.

\subsection{Stability and convergence}

\textbf{Lemma 1.} If
$g_{j}^{j+1}>g_{j-1}^{j+1}>\ldots>g_{0}^{j+1}>0$, $j=0,1,\ldots,M-1$
then for any function $v(t)$ defined on the mesh $\overline
\omega_{\tau}$ one has the inequalities
\begin{equation}\label{ur04}
v^{j+1}\left({_g}\Delta_{0t}^{\alpha}v\right)\geq
\frac{1}{2}{_g}\Delta_{0t}^{\alpha}(v^2)+\frac{1}{2g^{j+1}_j}\left({_g}\Delta_{0t}^{\alpha}v\right)^2,
\end{equation}
\begin{equation}\label{ur05}
v^{j}{_g}\Delta_{0t}^{\alpha}v\geq
\frac{1}{2}{_g}\Delta_{0t}^{\alpha}(v^2)-\frac{1}{2\left(g^{j+1}_j-g^{j+1}_{j-1}\right)}\left({_g}\Delta_{0t}^{\alpha}v\right)^2,
\end{equation}
where $g^{1}_{-1}=0$.

 {\bf Proof.} Let us consider the difference
$$
v^{j+1}{_g}\Delta_{0t}^{\alpha}v-
\frac{1}{2}{_g}\Delta_{0t}^{\alpha}(v^2)
$$
$$
=v^{j+1}\sum\limits_{s=0}^{j}
g_{s}^{j+1}(v^{s+1}-v^s)-\sum\limits_{s=0}^{j}
g_{s}^{j+1}(v^{s+1}-v^s)\left(\frac{v^{s+1}+v^{s}}{2}\right)
$$
$$
=\sum\limits_{s=0}^{j}
g_{s}^{j+1}(v^{s+1}-v^s)\left(v^{j+1}-\frac{v^{s+1}+v^{s}}{2}\right)
$$
$$
=\sum\limits_{s=0}^{j}
g_{s}^{j+1}(v^{s+1}-v^s)\left(\frac{1}{2}(v^{s+1}-v^s)+\sum\limits_{k=s+1}^{j}(v^{k+1}-v^k)\right)
$$
\begin{equation}\label{ur06}
=\frac{1}{2}\sum\limits_{s=0}^{j}
g_{s}^{j+1}(v^{s+1}-v^s)^2+\sum\limits_{k=1}^{j}(v^{k+1}-v^k)\sum\limits_{s=0}^{k-1}
 g_{s}^{j+1}(v^{s+1}-v^s).
\end{equation}
Here we consider the sums to be equal to zero if the upper summation
index is less than the lower one.

Let us introduce the following notation: $\sum_{s=0}^{k}
g_{s}^{j+1}(v^{s+1}-v^s)=w^{k+1}$, then $v^{1}-v^0=\left(
g_{0}^{j+1}\right)^{-1}w^{1}$, $v^{k+1}-v^k=\left(
g_{k}^{j+1}\right)^{-1}(w^{k+1}-w^{k})$, $k=1,2,\ldots,j$ and
rewrite the equality (\ref{ur06}) as
$$
\frac{1}{2}\left(g_{0}^{j+1}\right)^{-1}(w^{1})^2+
\frac{1}{2}\sum\limits_{k=1}^{j}\left(
g_{k}^{j+1}\right)^{-1}(w^{k+1}-w^{k})^2+\sum\limits_{k=1}^{j}\left(
g_{k}^{j+1}\right)^{-1}(w^{k+1}-w^{k})w^{k}
$$
$$=\frac{1}{2}\left(g_{0}^{j+1}\right)^{-1}(w^{1})^2+
\frac{1}{2}\sum\limits_{k=1}^{j}\left(
g_{k}^{j+1}\right)^{-1}\left((w^{k+1})^2-(w^{k})^2\right)
$$
$$=\frac{1}{2}\left( g_{j}^{j+1}\right)^{-1}(w^{j+1})^2+
\frac{1}{2}\sum\limits_{k=0}^{j-1}\frac{
g_{k+1}^{j+1}-g_{k}^{j+1}}{g_{k+1}^{j+1}
g_{k}^{j+1}}(w^{k+1})^2\geq\frac{1}{2}\left(
g_{j}^{j+1}\right)^{-1}(w^{j+1})^2,
$$
which is valid since $g_{k+1}^{j+1}-g_{k}^{j+1}>0$,
$k=0,1,\ldots,j-1$.

Let us prove now the inequality (\ref{ur05}). Since
$v^j=v^{j+1}-(v^{j+1}-v^{j})$, one obtains
$$
 v^{j}\Delta_{0t}^\alpha v - \frac{1}{2}\Delta_{0t}^\alpha (v^2)
 +\frac{1}{2\left(g_j^{j+1}-g_{j-1}^{j+1}\right)}(\Delta_{0t}^\alpha
 v)^2
$$
$$
=v^{j+1}\Delta_{0t}^\alpha v - \frac{1}{2}\Delta_{0t}^\alpha (v^2)
 +\frac{1}{2\left(g_j^{j+1}-g_{j-1}^{j+1}\right)}(\Delta_{0t}^\alpha
 v)^2-(v^{j+1}-v^{j})\Delta_{0t}^\alpha v
$$
$$
=\frac{1}{2}\left( g_{j}^{j+1}\right)^{-1}(w^{j+1})^2+
\frac{1}{2}\sum\limits_{k=0}^{j-1}\frac{
g_{k+1}^{j+1}-g_{k}^{j+1}}{g_{k+1}^{j+1} g_{k}^{j+1}}(w^{k+1})^2
$$
$$
+\frac{1}{2\left(g_j^{j+1}-g_{j-1}^{j+1}\right)}(w^{j+1})^2-\left(
g_{j}^{j+1}\right)^{-1}(w^{j+1}-w^{j})w^{j+1}
$$
$$
=\frac{g_{j-1}^{j+1}}{2g_{j}^{j+1}\left(g_j^{j+1}-g_{j-1}^{j+1}\right)}\left(w^{j+1}+
\frac{g_j^{j+1}-g_{j-1}^{j+1}}{g_{j-1}^{j+1}}w^{j}\right)^2+
\frac{1}{2}\sum\limits_{k=0}^{j-2}\frac{
g_{k+1}^{j+1}-g_{k}^{j+1}}{g_{k+1}^{j+1}
g_{k}^{j+1}}(w^{k+1})^2\geq0.
$$
The proof of the Lemma 1 is completed.

\textbf{Corollary 1.} If
$g_{j}^{j+1}>g_{j-1}^{j+1}>\ldots>g_{0}^{j+1}>0$ and
$\frac{g_{j}^{j+1}}{2g_{j}^{j+1}-g_{j-1}^{j+1}}\leq\sigma_{j+1}\leq1$,
where $j=0,1,\ldots,M-1$, $g_{-1}^{1}=0$, then for any function
$v(t)$ defined on the mesh $\overline\omega_{\tau}$ one has the
inequality
\begin{equation}\label{ur07}
 (\sigma_{j+1} v^{j+1}+(1-\sigma_{j+1})v^{j}){_g}\Delta_{0t}^\alpha v \geq \frac{1}{2}{_g}\Delta_{0t}^\alpha
 (v^2).
\end{equation}

\textbf{Theorem 1.} If
$$
g_{j}^{j+1}>g_{j-1}^{j+1}>\ldots>g_{0}^{j+1}\geq c_2>0, \quad
\frac{g_{j}^{j+1}}{2g_{j}^{j+1}-g_{j-1}^{j+1}}\leq\sigma_{j+1}\leq1,
$$
where  $j=0,1,\ldots,M-1$, $g_{-1}^{1}=0$, then the difference
scheme (\ref{ur03})--(\ref{ur03.1}) is unconditionally stable and
its solution satisfies the following a priori estimate:
\begin{equation}\label{ur08}
 \|y^{j+1}\|_0^2\leq\|y^0\|_0^2+\frac{1}{2\varkappa c_2}\max\limits_{0\leq j\leq
M}\|\varphi^{j}\|_0^2,
\end{equation}
where $(y,v)=\sum_{i=1}^{N-1}y_iv_ih$, $\|y\|_0^2=(y,y)$.

\textbf{Proof.} Taking the inner product of the equation
(\ref{ur03}) with $y^{(\sigma_{j+1})}$, we have
\begin{equation}\label{ur09}
 \left(y^{(\sigma_{j+1})},{_g}\Delta_{0t}^\alpha y\right)-\left(y^{(\sigma_{j+1})},\Lambda
 y^{(\sigma_{j+1})}\right)=\left(y^{(\sigma_{j+1})},\varphi^{j+1}\right).
\end{equation}

 Using inequality (\ref{ur07}) and the positive definiteness of operator $A=-\Lambda$ from  identity (\ref{ur09}) one obtains
\begin{equation}\label{ur010}
 \frac{1}{2}{_g}\Delta_{0t}^\alpha
 \|y\|_0^2+\varkappa\|y^{(\sigma_{j+1})}\|_0^2\leq
 \varepsilon\|y^{(\sigma_{j+1})}\|_0^2+\frac{1}{4\varepsilon}\|\varphi^{j+1}\|_0^2,\quad
 \varepsilon>0.
\end{equation}
From (\ref{ur010}), at $\varepsilon=\varkappa$ we get
\begin{equation}\label{ur011}
 {_g}\Delta_{0t}^\alpha
 \|y\|_0^2\leq
 \frac{1}{2\varkappa}\|\varphi^{j+1}\|_0^2.
\end{equation}
Let us rewrite inequality  (\ref{ur011}) in the form
\begin{equation}\label{ur011.11}
g_j^{j+1}\|y^{j+1}\|_0^2\leq\sum\limits_{s=1}^{j}\left(g_{s}^{j+1}-g_{s-1}^{j+1}\right)\|y^{s}\|_0^2+g_0^{j+1}\|y^0\|_0^2+\frac{1}{2\varkappa}\|\varphi^{j+1}\|_0^2.
\end{equation}
Noticing that $g_0^{j+1}\geq c_2>0$, we get
\begin{equation}\label{ur012}
g_j^{j+1}\|y^{j+1}\|_0^2\leq\sum\limits_{s=1}^{j}\left(g_{s}^{j+1}-g_{s-1}^{j+1}\right)\|y^{s}\|_0^2+g_0^{j+1}\left(\|y^0\|_0^2+\frac{1}{2\varkappa
c_2}\|\varphi^{j+1}\|_0^2\right).
\end{equation}
Denote
$$
E=\|y^0\|_0^2+\frac{1}{2\varkappa c_2}\max\limits_{0\leq j\leq
M}\|\varphi^{j}\|_0^2.
$$
The inequality (\ref{ur012}) is reduced to
\begin{equation}\label{ur013}
g_j^{j+1}\|y^{j+1}\|_0^2\leq\sum\limits_{s=1}^{j}\left(g_{s}^{j+1}-g_{s-1}^{j+1}\right)\|y^{s}\|_0^2+g_0^{j+1}E.
\end{equation}

It is obvious that at $j=0$ the a priori estimate (\ref{ur08})
follows from (\ref{ur013}). Let us prove that (\ref{ur08}) holds for
$j=1,2,\ldots$ by using the mathematical induction method. For this
purpose, let us assume that the a priori estimate (\ref{ur08}) takes
place for all $j=0,1,\ldots,k-1$:
$$
\|y^{j+1}\|_0^2\leq E, \quad j=0,1,\ldots, k-1.
$$

From (\ref{ur013}) at $j=k$ one has
\begin{equation}\label{ur014}
g_k^{k+1}\|y^{k+1}\|_0^2\leq\sum\limits_{s=1}^{k}\left(g_{s}^{k+1}-g_{s-1}^{k+1}\right)\|y^{s}\|_0^2+g_0^{k+1}E
$$
$$
\leq\sum\limits_{s=1}^{k}\left(g_{s}^{k+1}-g_{s-1}^{k+1}\right)E+g_0^{k+1}E=g_k^{k+1}E.
\end{equation}

The proof of Theorem 1 is completed.

A priori estimate  (\ref{ur08}) implies the stability of difference
scheme (\ref{ur03})--(\ref{ur03.1}).

\textbf{Theorem 2.} If the conditions of Theorem 1 are satisfied and
difference scheme (\ref{ur03})--(\ref{ur03.1}) has the approximation
order  $\mathcal{O}(N^{-r_1}+M^{-r_2})$, where $r_1$ and $r_2$ are
some known positive numbers, then the solution of difference scheme
(\ref{ur03})--(\ref{ur03.1}) converges to the solution of
differential problem (\ref{ur01})--(\ref{ur02}) in the mesh  $L_2$ -
norm with the rate equal to the order of the approximation error
$\mathcal{O}(N^{-r_1}+M^{-r_2})$.

\textbf{Proof.}  Let us introduce the error  $z=y-u$ and substitute it
into (\ref{ur03})--(\ref{ur03.1}). Then we obtain the problem for
the error
\begin{equation}\label{ur015}
{_g}\Delta_{0t}^{\alpha}z_i=\Lambda z^{(\sigma_{j+1})}_i
+\psi_i^{j+1}, \quad i=1,\ldots,N-1,\quad j=0,1,\ldots,M-1,
\end{equation}
\begin{equation}
z(0,t)=0,\quad z(l,t)=0,\quad t\in \overline \omega_{\tau}, \quad
z(x,0)=0,\quad  x\in \overline \omega_{h},\label{ur016}
\end{equation}
where $\psi_i^{j+1}=\Lambda
u^{(\sigma_{j+1})}_i-{_g}\Delta_{0t}^{\alpha}u_i+\varphi_i^{j+1}$,\quad
$\psi_i^{j+1}=\mathcal{O}(N^{-r_1}+M^{-r_2})$.

Since the conditions of Theorem 1 are fulfilled, then a priori estimate
(\ref{ur08}) holds true for the solution of problem
(\ref{ur015})--(\ref{ur016}) and, therefore, the following
inequality takes place
$$
\|z\|_0\leq\frac{1}{\sqrt{2\varkappa c_2}}\max\limits_{0\leq j\leq
M}\|\psi^{j}\|_0=\mathcal{O}(N^{-r_1}+M^{-r_2}),
$$
which implies the convergence in the mesh $L_2$ - norm with the rate
$\mathcal{O}(N^{-r_1}+M^{-r_2})$.

%For example, using the discrete Green's formula and noticing that
%$y(0,t)=0$, $y(l,t)=0$, we have

\section{A new $\bf L2-1_\sigma$ fractional numerical differentiation formula}

In this section a difference analog of the Caputo fractional
derivative with the  approximation order $O(\tau^{3-\alpha})$ is
constructed and its basic properties are investigated.

Let us consider the uniform mesh  $\bar \omega_{\tau}=\{t_j=j\tau,\,
j=0,1,\ldots,M; \,T=\tau M\}$. Let $\sigma=1-\frac{\alpha}{2}$, then
for the Caputo fractional derivative of the order  $\alpha$,
$0<\alpha<1$, of the function  $u(t)\in \mathcal{C}^{3}[0,T]$  at
the fixed point  $t_{j+\sigma}$,  $j\in \{0,1,\ldots,M-1\}$ the
following equalities hold
$$
\partial_{0t_{j+\sigma}}^{\alpha}u(\eta)=\frac{1}{\Gamma(1-\alpha)}\int\limits_{0}^{t_{j+\sigma}}\frac{u'(\eta)d\eta}{(t_{j+\sigma}-\eta)^{\alpha}}
$$
\begin{equation}\label{ur0.99}
=
\frac{1}{\Gamma(1-\alpha)}\sum\limits_{s=1}^{j}\int\limits_{t_{s-1}}^{t_{s}}\frac{u'(\eta)d\eta}{(t_{j+\sigma}-\eta)^{\alpha}}+
\frac{1}{\Gamma(1-\alpha)}\int\limits_{t_{j}}^{t_{j+\sigma}}\frac{u'(\eta)d\eta}{(t_{j+\sigma}-\eta)^{\alpha}}.
\end{equation}

As in \cite{Sun6}, on each interval $[t_{s-1},t_s]$ ($1\leq
s\leq j$), denoting the quadratic interpolation ${\Pi}_{2,s}u(t)$ of
$u(t)$ using three points $(t_{s-1},u(t_{s-1}))$, $(t_{s},u(t_{s}))$
and $(t_{s+1},u(t_{s+1}))$, we get
$$
{\Pi}_{2,s}u(t)=u(t_{s-1})\frac{(t-t_{s})(t-t_{s+1})}{2\tau^2}
$$
$$
-u(t_{s})\frac{(t-t_{s-1})(t-t_{s+1})}{\tau^2}+u(t_{s+1})\frac{(t-t_{s-1})(t-t_{s})}{2\tau^2},
$$
\begin{equation}\label{ur0.991}
\left({\Pi}_{2,s}u(t)\right)'=u_{t,s}+u_{\bar
tt,s}(t-t_{s+1/2})=u_{t,s-1}+u_{\bar tt,s}(t-t_{s-1/2}),
\end{equation}
and
\begin{equation}\label{ur0.992}
u(t)-{\Pi}_{2,s}u(t)=\frac{u'''(\bar\xi_s)}{6}(t-t_{s-1})(t-t_{s})(t-t_{s+1}),
\end{equation}
where $t\in[t_{s-1},t_{s+1}]$, $\bar\xi_s\in(t_{s-1},t_{s+1})$,
$u_{t,s}=(u(t_{s+1})-u(t_s))/\tau$, $u_{\bar
t,s}=(u(t_{s})-u(t_{s-1}))/\tau$.

In (\ref{ur0.99}), we use ${\Pi}_{2,s}u(t)$ to approximate $u(t)$ on
the interval $[t_{s-1},t_{s}]$ ($1\leq s\leq j$). Taking into
account the equality
\begin{equation}\label{ur0.993}
\int\limits_{t_{s-1}}^{t_s}(\eta-t_{s-1/2})(t_{j+\sigma}-\eta)^{-\alpha}d\eta=\frac{\tau^{2-\alpha}}{1-\alpha}b_{j-s+1}^{(\alpha,\sigma)},
\quad 1\leq s\leq j
\end{equation}
with
$$
b_{l}^{(\alpha,\sigma)}=\frac{1}{2-\alpha}\left[(l+\sigma)^{2-\alpha}-(l-1+\sigma)^{2-\alpha}\right]-\frac{1}{2}\left[(l+\sigma)^{1-\alpha}+(l-1+\sigma)^{1-\alpha}\right],
$$
$l\geq 1$, from (\ref{ur0.99}) and  (\ref{ur0.991}) we obtain the
difference analog of the Caputo fractional derivative of the order
$\alpha$ ($0<\alpha<1$) for the function $u(t)$ in the following
form:
$$
\partial_{0t_{j+\sigma}}^{\alpha}u(\eta)=
\frac{1}{\Gamma(1-\alpha)}\sum\limits_{s=1}^{j}\int\limits_{t_{s-1}}^{t_{s}}\frac{u'(\eta)d\eta}{(t_{j+\sigma}-\eta)^{\alpha}}+
\frac{1}{\Gamma(1-\alpha)}\int\limits_{t_{j}}^{t_{j+\sigma}}\frac{u'(\eta)d\eta}{(t_{j+\sigma}-\eta)^{\alpha}}
$$
$$
\approx
\frac{1}{\Gamma(1-\alpha)}\sum\limits_{s=1}^{j}\int\limits_{t_{s-1}}^{t_{s}}\frac{\left({\Pi}_{2,s}u(\eta)\right)'d\eta}{(t_{j+\sigma}-\eta)^{\alpha}}+
\frac{u_{t,j}}{\Gamma(1-\alpha)}\int\limits_{t_{j}}^{t_{j+\sigma}}\frac{d\eta}{(t_{j+\sigma}-\eta)^{\alpha}}
$$
$$
=\frac{1}{\Gamma(1-\alpha)}\sum\limits_{s=1}^{j}\int\limits_{t_{s-1}}^{t_{s}}\frac{u_{t,s-1}+u_{\bar
tt,s}(\eta-t_{s-1/2})d\eta}{(t_{j+\sigma}-\eta)^{\alpha}}+
\frac{u_{t,j}}{\Gamma(1-\alpha)}\int\limits_{t_{j}}^{t_{j+\sigma}}\frac{d\eta}{(t_{j+\sigma}-\eta)^{\alpha}}
$$
$$
=\frac{\tau^{1-\alpha}}{\Gamma{(2-\alpha)}}\left(\sum\limits_{s=1}^{j}\left(a_{j-s+1}^{(\alpha,\sigma)}u_{t,s-1}+b_{j-s+1}^{(\alpha,\sigma)}u_{\bar
tt,s}\tau\right)+a_0^{(\alpha,\sigma)}u_{t,j}\right)
$$
$$
=\frac{\tau^{1-\alpha}}{\Gamma{(2-\alpha)}}\left(\sum\limits_{s=1}^{j}\left(a_{j-s+1}^{(\alpha,\sigma)}u_{t,s-1}+b_{j-s+1}^{(\alpha,\sigma)}(u_{t,s}-u_{t,s-1})\right)+a_0^{(\alpha,\sigma)}u_{t,j}\right)
$$
\begin{equation}
=\frac{\tau^{1-\alpha}}{\Gamma{(2-\alpha)}}\sum\limits_{s=0}^{j}c_{j-s}^{(\alpha,\sigma)}u_{t,s}=\Delta_{0t_{j+\sigma}}^\alpha
u, \label{ur0.994}
\end{equation}
where
$$
a_{0}^{(\alpha,\sigma)}=\sigma^{1-\alpha},\quad
a_{l}^{(\alpha,\sigma)}=(l+\sigma)^{1-\alpha}-(l-1+\sigma)^{1-\alpha},
\quad l\geq 1;
$$
$c_{0}^{(\alpha,\sigma)}=a_{0}^{(\alpha,\sigma)}$, for $j=0$; and
for $j\geq 1$,
\begin{equation}
c_{s}^{(\alpha,\sigma)}=
\begin{cases}
a_0^{(\alpha,\sigma)}+b_{1}^{(\alpha,\sigma)}, \quad\quad\quad \quad\,\,\, s=0,\\
a_{s}^{(\alpha,\sigma)}+b_{s+1}^{(\alpha,\sigma)}-b_{s}^{(\alpha,\sigma)}, \quad 1\leq s\leq j-1,\\
a_{j}^{(\alpha,\sigma)}-b_{j}^{(\alpha,\sigma)},
\quad\quad\quad\quad\,\,\, s=j. \label{ur0.995}
\end{cases}
\end{equation}

We call the fractional numerical differentiation formula
(\ref{ur0.994}) for the Caputo fractional derivative of order
$\alpha$ ($0<\alpha<1$) the $L2$-$1_\sigma$ formula.

%Now, truncation errors of the new $L1$-$2_\sigma$ formula
%(\ref{ur0.994}) are

\textbf{Lemma 2.} For any $\alpha\in(0,1)$ and $u(t)\in
\mathcal{C}^3[0,t_{j+1}]$
\begin{equation}
|\partial_{0t_{j+\sigma}}^{\alpha}u-\Delta_{0t_{j+\sigma}}^\alpha
u|=\mathcal{O}(\tau^{3-\alpha}).
 \label{ur0.996}
\end{equation}
\textbf{Proof.} Let
$\partial_{0t_{j+\sigma}}^{\alpha}u-\Delta_{0t_{j+\sigma}}^\alpha
u=R_{1}^{j}+R_{j}^{j+\sigma}$, where
$$
R_{1}^{j}=\frac{1}{\Gamma(1-\alpha)}\sum\limits_{s=1}^{j}\int\limits_{t_{s-1}}^{t_{s}}\frac{u'(\eta)d\eta}{(t_{j+\sigma}-\eta)^{\alpha}}-
\frac{1}{\Gamma(1-\alpha)}\sum\limits_{s=1}^{j}\int\limits_{t_{s-1}}^{t_{s}}\frac{\left({\Pi}_{2,s}u(\eta)\right)'d\eta}{(t_{j+\sigma}-\eta)^{\alpha}}
$$
$$
=\frac{1}{\Gamma(1-\alpha)}\sum\limits_{s=1}^{j}\int\limits_{t_{s-1}}^{t_{s}}\left(u(\eta)-{\Pi}_{2,s}u(\eta)\right)'{(t_{j+\sigma}-\eta)^{-\alpha}}d\eta
$$
$$
=-\frac{\alpha}{\Gamma(1-\alpha)}\sum\limits_{s=1}^{j}\int\limits_{t_{s-1}}^{t_{s}}\left(u(\eta)-{\Pi}_{2,s}u(\eta)\right){(t_{j+\sigma}-\eta)^{-\alpha-1}}d\eta
$$
$$
=-\frac{\alpha}{6\Gamma(1-\alpha)}\sum\limits_{s=1}^{j}\int\limits_{t_{s-1}}^{t_{s}}u'''(\bar\xi_s)(\eta-t_{s-1})(\eta-t_{s})(\eta-t_{s+1}){(t_{j+\sigma}-\eta)^{-\alpha-1}}d\eta,
$$
$$
R_{j}^{j+\sigma}=\frac{1}{\Gamma(1-\alpha)}\int\limits_{t_{j}}^{t_{j+\sigma}}\frac{u'(\eta)d\eta}{(t_{j+\sigma}-\eta)^{\alpha}}-
\frac{u_{t,j}}{\Gamma(1-\alpha)}\int\limits_{t_{j}}^{t_{j+\sigma}}\frac{d\eta}{(t_{j+\sigma}-\eta)^{\alpha}}
$$
$$
=\frac{1}{\Gamma(1-\alpha)}\int\limits_{t_{j}}^{t_{j+\sigma}}\frac{(u'(\eta)-u_{t,j})d\eta}{(t_{j+\sigma}-\eta)^{\alpha}}=
\frac{1}{\Gamma(1-\alpha)}\int\limits_{t_{j}}^{t_{j+\sigma}}\frac{(u'(t_{j+1/2})-u_{t,j})d\eta}{(t_{j+\sigma}-\eta)^{\alpha}}+
$$
$$
+\frac{u''(t_{j+1/2})}{\Gamma(1-\alpha)}\int\limits_{t_{j}}^{t_{j+\sigma}}\frac{(\eta-t_{j+1/2})d\eta}{(t_{j+\sigma}-\eta)^{\alpha}}+\mathcal{O}(\tau^{3-\alpha})
$$
$$
=\frac{u''(t_{j+1/2})}{\Gamma(1-\alpha)}\int\limits_{t_{j}}^{t_{j+\sigma}}\frac{(\eta-t_{j+1/2})d\eta}{(t_{j+\sigma}-\eta)^{\alpha}}+\mathcal{O}(\tau^{3-\alpha}).
$$
We estimate the error  $R_{1}^{j}$ similarly to  \cite{Sun6}:
$$
|R_{1}^{j}|\leq\frac{\alpha|u'''(\xi)|}{6\Gamma(1-\alpha)}\sum\limits_{s=1}^{j}\int\limits_{t_{s-1}}^{t_{s}}(\eta-t_{s-1})(t_{s}-\eta)(t_{s+1}-\eta){(t_{j+\sigma}-\eta)^{-\alpha-1}}d\eta
$$
$$
\leq\frac{\alpha|u'''(\xi)|\tau^3}{3\Gamma(1-\alpha)}\sum\limits_{s=1}^{j}\int\limits_{t_{s-1}}^{t_{s}}{(t_{j+\sigma}-\eta)^{-\alpha-1}}d\eta=
\frac{\alpha|u'''(\xi)|\tau^3}{3\Gamma(1-\alpha)}\int\limits_{0}^{t_{j}}{(t_{j+\sigma}-\eta)^{-\alpha-1}}d\eta
$$
$$
=\frac{|u'''(\xi)|\tau^3}{3\Gamma(1-\alpha)}\left(\frac{1}{\sigma^\alpha\tau^\alpha}-\frac{1}{(j+\sigma)^\alpha\tau^\alpha}\right)
\leq\frac{|u'''(\xi)|}{3\sigma^\alpha\Gamma(1-\alpha)}\tau^{3-\alpha},
\quad \xi\in(0,t_j).
$$

Since
$$
\int\limits_{t_{j}}^{t_{j+\sigma}}\frac{(\eta-t_{j+1/2})d\eta}{(t_{j+\sigma}-\eta)^{\alpha}}=
\frac{\tau
t_\sigma^{1-\alpha}\left(2\sigma+\alpha-2\right)}{2(1-\alpha)(2-\alpha)}=0
$$
the error  $|R_{j}^{j+\sigma}|=\mathcal{O}(\tau^{3-\alpha})$. Lemma~2 is proved.

\subsection{Basic properties of the new $\bf L2-1_\sigma$ fractional numerical differentiation formula.}

\textbf{Lemma 3.} For all  $s=1,2,\ldots$ and $0<\alpha<1$ the
following inequalities hold
$$
\frac{1}{2}<\varkappa_s<\frac{1}{2-\alpha},
$$
where
$$
\varkappa_s=\frac{(s+\sigma)^{2-\alpha}-(s-1+\sigma)^{2-\alpha}-(2-\alpha)(s-1+\sigma)^{1-\alpha}}{(2-\alpha)((s+\sigma)^{1-\alpha}-(s-1+\sigma)^{1-\alpha})}.
$$

\textbf{Proof.} Let us consider two functions
$$
f_\alpha(x)=\frac{(x+1)^{2-\alpha}-x^{2-\alpha}-(2-\alpha)x^{1-\alpha}}{(2-\alpha)((x+1)^{1-\alpha}-x^{1-\alpha})}=
\int_{0}^{1}\frac{(z+x)^{1-\alpha}-x^{1-\alpha}}{(1+x)^{1-\alpha}-x^{1-\alpha}}dz,\quad
x>0
$$
 and
$$
g_\alpha(z,x)=\frac{(z+x)^{1-\alpha}-x^{1-\alpha}}{(1+x)^{1-\alpha}-x^{1-\alpha}}=\frac{z\int\limits_{0}^{1}\frac{d\xi}{(x+z\xi)^{\alpha}}}{\int\limits_{0}^{1}\frac{d\xi}{(x+\xi)^{\alpha}}},
\quad 0<z<1,\quad x>0.
$$
For all  $x>0$ and $0<z<1$ the following inequalities hold
$$
\int\limits_{0}^{1}\frac{d\xi}{(x+\xi)^{\alpha}}<\int\limits_{0}^{1}\frac{d\xi}{(x+z\xi)^{\alpha}}<\int\limits_{0}^{1}\frac{d\xi}{(zx+z\xi)^{\alpha}}=z^{-\alpha}\int\limits_{0}^{1}\frac{d\xi}{(x+\xi)^{\alpha}}.
$$
Therefore, for the function  $g_\alpha(z,x)$ for all  $x>0$ and
$0<z<1$ the inequalities
\begin{equation}\label{ur4.101}
z<g_\alpha(z,x)<z^{1-\alpha}
\end{equation}
are valid.

 Integrating  (\ref{ur4.101}) with respect to
$z$ from  $0$ to $1$, we get the inequalities
$$
\frac{1}{2}<f_\alpha(x)<\frac{1}{2-\alpha},
$$
which hold for all  $x>0$. Lemma 3 is proved.

\textbf{Corollary 2.} For any $\alpha$ ($0<\alpha<1$), it holds
$b_{s}^{(\alpha,\sigma)}>0$, $s\geq1$.

The latter follows from the equality
$$
b_{s}^{(\alpha,\sigma)}=\left[(s+\sigma)^{1-\alpha}-(s-1+\sigma)^{1-\alpha}\right]\left(\varkappa_s-\frac{1}{2}\right).
$$

\textbf{Lemma 4.} For any $\alpha$ ($0<\alpha<1$) and
$c_s^{(\alpha,\sigma)}$ ($0\leq s\leq j, j\geq1$) defined in
(\ref{ur0.995}), it holds
\begin{equation}\label{ur4.102}
c_{j}^{(\alpha,\sigma)}>\frac{1-\alpha}{2}(j+\sigma)^{-\alpha},
\end{equation}
\begin{equation}\label{ur4.103}
c_{0}^{(\alpha,\sigma)}>c_{1}^{(\alpha,\sigma)}>c_{2}^{(\alpha,\sigma)}>\ldots>c_{j-1}^{(\alpha,\sigma)}>c_{j}^{(\alpha,\sigma)},
\end{equation}
\begin{equation}\label{ur4.104}
(2\sigma-1)c_{0}^{(\alpha,\sigma)}-\sigma c_{1}^{(\alpha,\sigma)}>0,
\end{equation}
where $\sigma=1-\alpha/2$.

 \textbf{Proof.}
For $j\geq 1$ we get
$$
c_{j}^{(\alpha,\sigma)}=
a_{j}^{(\alpha,\sigma)}-b_{j}^{(\alpha,\sigma)}=\left((j+\sigma)^{1-\alpha}-(j-1+\sigma)^{1-\alpha}\right)\left(\frac{3}{2}-\varkappa_j\right)
$$
$$
>\left((j+\sigma)^{1-\alpha}-(j-1+\sigma)^{1-\alpha}\right)\left(\frac{3}{2}-\frac{1}{2-\alpha}\right)
$$
$$
>\frac{1-\alpha}{2}\int\limits_{0}^{1}\frac{d\eta}{(j+\sigma-\eta)^\alpha}>\frac{1-\alpha}{2}(j+\sigma)^{-\alpha}.
$$
Inequality (\ref{ur4.102}) is proved. Let us prove
inequality (\ref{ur4.103}).

For $1\leq s\leq j-2$ ($j\geq 3$) we have
$$
c_{s}^{(\alpha,\sigma)}-c_{s+1}^{(\alpha,\sigma)}=a_{s}^{(\alpha,\sigma)}-a_{s+1}^{(\alpha,\sigma)}+
2b_{s+1}^{(\alpha,\sigma)}-b_{s}^{(\alpha,\sigma)}-b_{s+2}^{(\alpha,\sigma)}
$$
$$
=\frac{1}{2}\left((s+2+\sigma)^{1-\alpha}-3(s+1+\sigma)^{1-\alpha}+3(s+\sigma)^{1-\alpha}-(s-1+\sigma)^{1-\alpha}\right)
$$
$$
+\frac{1}{2-\alpha}\left(-(s+2+\sigma)^{2-\alpha}+3(s+1+\sigma)^{2-\alpha}-3(s+\sigma)^{2-\alpha}+(s-1+\sigma)^{2-\alpha}\right)
$$
$$
=\frac{\alpha(1-\alpha)(1+\alpha)}{2}\int\limits_{0}^{1}dz_1\int\limits_{0}^{1}dz_2\int\limits_{0}^{1}\frac{dz_3}{(s-1+\sigma+z_1+z_2+z_3)^{\alpha+2}}
$$
$$
+\alpha(1-\alpha)\int\limits_{0}^{1}dz_1\int\limits_{0}^{1}dz_2\int\limits_{0}^{1}\frac{dz_3}{(s-1+\sigma+z_1+z_2+z_3)^{\alpha+1}}
$$
$$
>\frac{\alpha(1-\alpha)(1+\alpha)}{2}(s+2+\sigma)^{-\alpha-2}+\alpha(1-\alpha)(s+2+\sigma)^{-\alpha-1}>0.
$$

For $s=j-1$ ($j\geq 2$) we get
$$
c_{s}^{(\alpha,\sigma)}-c_{s+1}^{(\alpha,\sigma)}=c_{j-1}^{(\alpha,\sigma)}-c_{j}^{(\alpha,\sigma)}=a_{j-1}^{(\alpha,\sigma)}-a_{j}^{(\alpha,\sigma)}+
2b_{j}^{(\alpha,\sigma)}-b_{j-1}^{(\alpha,\sigma)}
$$
$$
>a_{j-1}^{(\alpha,\sigma)}-a_{j}^{(\alpha,\sigma)}+
2b_{j}^{(\alpha,\sigma)}-b_{j-1}^{(\alpha,\sigma)}-b_{j+1}^{(\alpha,\sigma)}
$$
$$
>\frac{\alpha(1-\alpha)(1+\alpha)}{2}(j+1+\sigma)^{-\alpha-2}+\alpha(1-\alpha)(j+1+\sigma)^{-\alpha-1}>0.
$$

For inequality  (\ref{ur4.103}) it remains to prove the case  $s=0$,
that is  $c_{0}^{(\alpha,\sigma)}>c_{1}^{(\alpha,\sigma)}$ which obviously follows from  (\ref{ur4.104}).
It is enough to prove inequality  (\ref{ur4.104}).

For $j=1$ we get
$$
(2\sigma-1)c_{0}^{(\alpha,\sigma)}-\sigma c_{1}^{(\alpha,\sigma)}=
(2\sigma-1)(a_{0}^{(\alpha,\sigma)}+b_{1}^{(\alpha,\sigma)})-\sigma
(a_{1}^{(\alpha,\sigma)}-b_{1}^{(\alpha,\sigma)})
$$
$$
=\left(\frac{2\sigma-1}{2\sigma}-\frac{2\sigma-1}{2}\right)(1+\sigma)^{1-\alpha}=\frac{(2\sigma-1)(1-\sigma)}{2\sigma}(1+\sigma)^{1-\alpha}>0.
$$

For $j\geq2$ we get
$$
(2\sigma-1)c_{0}^{(\alpha,\sigma)}-\sigma c_{1}^{(\alpha,\sigma)}=
(2\sigma-1)(a_{0}^{(\alpha,\sigma)}+b_{1}^{(\alpha,\sigma)})-\sigma
(a_{1}^{(\alpha,\sigma)}+b_{2}^{(\alpha,\sigma)}-b_{1}^{(\alpha,\sigma)})
$$
$$
=\frac{4\sigma-1}{2\sigma}(1+\sigma)^{1-\alpha}-(2+\sigma)^{1-\alpha}=
(1+\sigma)^{1-\alpha}\left(\frac{4\sigma-1}{2\sigma}-\left(1+\frac{1}{1+\sigma}\right)^{1-\alpha}\right)
$$
$$
>(1+\sigma)^{1-\alpha}\left(\frac{4\sigma-1}{2\sigma}-1-\frac{1-\alpha}{1+\sigma}\right)=\frac{(2\sigma-1)(1-\sigma)}{2\sigma(1+\sigma)^\alpha}>0.
$$
Here we used the inequality  $(1+t)^\gamma<1+\gamma t$ which is
valid for all  $t>0$ and $0<\gamma<1$. Lemma 4 is proved.

\subsection{Test example}

In this subsection, the validity and numerical accuracy of the new
presented  $L2$-$1_\sigma$ formula (\ref{ur0.994}) are demonstrated
by a test example.

Let us take a positive integer $M$, let  $\tau=1/(M-1+\sigma)$ and
denote
$$
E_{L2-1_\sigma}^{M}(\tau)=|\partial_{0t_{M-1+\sigma}}^{\alpha}f(t)-\Delta_{0t_{M-1+\sigma}}^{\alpha}f(t)|.
$$
\textbf{Example.} Let $f(t)=t^{4+\alpha}, \quad 0<\alpha<1$. Compute
the $\alpha$-order Caputo fractional derivative of $f(t)$ at
$t=t_{M-1+\sigma}=1$ numerically.

The exact solution is given by
$$
\left.\partial_{0t}^{\alpha}t^{4+\alpha}\right|_{t=1}=\frac{\Gamma(5+\alpha)}{24}.
$$

Taking different temporal stepsizes  $M=10$, \, $20$, \, $40$, \,
$80$, \, $160$, \, $320$, \, $640$, \, $1280$, \, $2560$, \, $5120$,
\, we compute the example using  $L2-1_\sigma$ formula
(\ref{ur0.994}) and compare the results with those obtained with the
help of the  $L1-2$ formula in  {\cite{Sun6}}. Table 1 lists the
computational errors and numerical convergence order ($CO$) at
$t_{M-1+\sigma}=1$ with different parameters $\alpha=0.9$,\,
$0.5$,\, $0.1$.

\begin{tabular}{lc}
{\bf Table 1.}\\
Computational errors and convergence order with different
\\ temporal stepsizes \\
 \hline
$\alpha$ \hspace{6mm} $M$  \hspace{10mm} $E_{L1-2}^{M}(\tau)$\cite{Sun6} \hspace{8mm} $CO_{E_{L1-2}^{M}}$    \hspace{8mm}{$E_{L2-1_\sigma}^{M}(\tau)$} \hspace{10mm} $CO_{E_{L2-1_\sigma}^{M}}$ \\
\hline
0.9       \hspace{3mm}  10        \hspace{10mm}  $1.070471e-1$   \hspace{15mm}           \hspace{13mm}   $1.922978e-2$                                    \\
          \hspace{10mm} 20        \hspace{10mm}  $2.699702e-2$   \hspace{8mm}    $1.99$  \hspace{11mm}   $4.368964e-3$      \hspace{8mm} $2.07$         \\
          \hspace{10mm} 40        \hspace{10mm}  $6.545547e-3$   \hspace{8mm}    $2.04$  \hspace{11mm}   $1.009364e-3$      \hspace{8mm} $2.08$         \\
          \hspace{10mm} 80        \hspace{10mm}  $1.556707e-3$   \hspace{8mm}    $2.07$  \hspace{11mm}   $2.347614e-4$      \hspace{8mm} $2.09$         \\
          \hspace{10mm} 160       \hspace{8mm}   $3.666902e-4$   \hspace{8mm}    $2.09$  \hspace{11mm}   $5.473732e-5$      \hspace{8mm} $2.09$         \\
          \hspace{10mm} 320       \hspace{8mm}   $8.595963e-5$   \hspace{8mm}    $2.09$  \hspace{11mm}   $1.277246e-5$      \hspace{8mm} $2.10$         \\
          \hspace{10mm} 640       \hspace{8mm}   $2.010152e-5$   \hspace{8mm}    $2.10$  \hspace{11mm}   $2.980723e-6$      \hspace{8mm} $2.10$         \\
          \hspace{10mm} 1280      \hspace{6mm}   $4.694884e-6$   \hspace{8mm}    $2.10$  \hspace{11mm}   $6.955612e-7$      \hspace{8mm} $2.10$         \\
          \hspace{10mm} 2560      \hspace{6mm}   $1.095840e-6$   \hspace{8mm}    $2.10$  \hspace{11mm}   $1.622925e-7$      \hspace{8mm} $2.10$         \\
          \hspace{10mm} 5120      \hspace{6mm}   $2.556990e-7$   \hspace{8mm}    $2.10$  \hspace{11mm}   $3.786340e-8$      \hspace{8mm} $2.10$         \\
         \\
0.5       \hspace{3mm}  10        \hspace{10mm}  $1.350657e-2$   \hspace{15mm}           \hspace{13mm}   $3.756950e-3$                                    \\
          \hspace{10mm} 20        \hspace{10mm}  $2.612085e-3$   \hspace{8mm}    $2.37$  \hspace{11mm}   $7.231988e-4$      \hspace{8mm} $2.33$         \\
          \hspace{10mm} 40        \hspace{10mm}  $4.861786e-4$   \hspace{8mm}    $2.43$  \hspace{11mm}   $1.367574e-4$      \hspace{8mm} $2.38$         \\
          \hspace{10mm} 80        \hspace{10mm}  $8.864502e-5$   \hspace{8mm}    $2.46$  \hspace{11mm}   $2.544814e-5$      \hspace{8mm} $2.42$         \\
          \hspace{10mm} 160       \hspace{8mm}   $1.597499e-5$   \hspace{8mm}    $2.47$  \hspace{11mm}   $4.673501e-6$      \hspace{8mm} $2.44$         \\
          \hspace{10mm} 320       \hspace{8mm}   $2.859085e-6$   \hspace{8mm}    $2.48$  \hspace{11mm}   $8.495470e-7$      \hspace{8mm} $2.46$         \\
          \hspace{10mm} 640       \hspace{8mm}   $5.095342e-7$   \hspace{8mm}    $2.49$  \hspace{11mm}   $1.532461e-7$      \hspace{8mm} $2.47$         \\
          \hspace{10mm} 1280      \hspace{6mm}   $9.056389e-8$   \hspace{8mm}    $2.49$  \hspace{11mm}   $2.748687e-8$      \hspace{8mm} $2.48$         \\
          \hspace{10mm} 2560      \hspace{6mm}   $1.606869e-8$   \hspace{8mm}    $2.49$  \hspace{11mm}   $4.909831e-9$      \hspace{8mm} $2.48$         \\
          \hspace{10mm} 5120      \hspace{6mm}   $2.847764e-9$   \hspace{8mm}    $2.50$  \hspace{11mm}   $8.743961e-10$     \hspace{6mm} $2.49$         \\
         \\
0.1       \hspace{3mm}  10        \hspace{10mm}  $6.238229e-4$   \hspace{15mm}           \hspace{13mm}   $2.686107e-4$                                    \\
          \hspace{10mm} 20        \hspace{10mm}  $9.663202e-5$   \hspace{8mm}    $2.69$  \hspace{11mm}   $4.492624e-5$      \hspace{8mm} $2.57$         \\
          \hspace{10mm} 40        \hspace{10mm}  $1.444281e-5$   \hspace{8mm}    $2.74$  \hspace{11mm}   $7.204745e-6$      \hspace{8mm} $2.64$         \\
          \hspace{10mm} 80        \hspace{10mm}  $2.111896e-6$   \hspace{8mm}    $2.77$  \hspace{11mm}   $1.119177e-6$      \hspace{8mm} $2.68$         \\
          \hspace{10mm} 160       \hspace{8mm}   $3.043133e-7$   \hspace{8mm}    $2.79$  \hspace{11mm}   $1.696376e-7$      \hspace{8mm} $2.72$         \\
          \hspace{10mm} 320       \hspace{8mm}   $4.338827e-8$   \hspace{8mm}    $2.81$  \hspace{11mm}   $2.522442e-8$      \hspace{8mm} $2.75$         \\
          \hspace{10mm} 640       \hspace{8mm}   $6.136347e-9$   \hspace{8mm}    $2.82$  \hspace{11mm}   $3.694254e-9$      \hspace{8mm} $2.77$         \\
          \hspace{10mm} 1280      \hspace{6mm}   $8.622698e-10$  \hspace{6mm}    $2.83$  \hspace{11mm}   $5.344856e-10$     \hspace{6mm} $2.79$         \\
          \hspace{10mm} 2560      \hspace{6mm}   $1.205229e-10$  \hspace{6mm}    $2.84$  \hspace{11mm}   $7.656497e-11$     \hspace{6mm} $2.80$         \\
          \hspace{10mm} 5120      \hspace{6mm}   $1.676992e-11$  \hspace{6mm}    $2.85$  \hspace{11mm}   $1.087796e-11$     \hspace{6mm} $2.82$         \\
\hline
\end{tabular}

\section{A second order difference scheme for the time fractional diffusion
equation}

%A difference scheme for the time fractional diffusion equation of
%the second order of approximation. Stability and convergence.
%Numerical results}

In this section for problem  (\ref{ur01})--(\ref{ur02})   a
difference scheme with the approximation order
$\mathcal{O}(h^2+\tau^2)$ is constructed. The stability of the
constructed difference scheme as well as its convergence in the mesh
$L_2$ - norm with the rate equal to the order of the approximation
error is proved. The obtained results are supported with numerical
calculations carried out for a test example.

\subsection{Derivation  of the difference scheme}

\textbf{Lemma 5.} For any functions  $k_1(x)\in \mathcal{C}_{x}^{3}$
and $v(x)\in\mathcal{C}_{x}^{4}$ the following equality is valid:
$$
\left.\frac{d}{dx}\left(k_1(x)\frac{d}{dx}v(x)\right)\right|_{x=x_i}
$$
\begin{equation}\label{ur16.1}
=\frac{k_1(x_{i+1/2})v(x_{i+1})-(k_1(x_{i+1/2})+k_1(x_{i-1/2}))v(x_{i})+k_1(x_{i-1/2})v(x_{i-1})}{h^2}+\mathcal{O}(h^2).
\end{equation}

%\textbf{Proof.} Using the Taylor series expansions of the
%functions we obtain the following equalities:
%\begin{equation}\label{ur16.2}
%\frac{v(x_{i+1})-v(x_i)}{h}=v'(x_i)+v''(x_i)\frac{h}{2}+v'''(x_i)\frac{h^2}{6}+\mathcal{O}(h^3).
%\end{equation}
%\begin{equation}\label{ur16.3}
%k_1(x_{i+1/2})=k_1(x_i)+k'_1(x_i)\frac{h}{2}+k''_1(x_i)\frac{h^2}{8}+\mathcal{O}(h^3).
%\end{equation}
%From (\ref{ur16.2}) and (\ref{ur16.3}) we get
%$$
%k_1(x_{i+1/2})\frac{v(x_{i+1})-v(x_i)}{h}=k_1(x_i)v'(x_i)+\frac{h}{2}\left(k_1v'\right)'(x_i)+
%$$
%\begin{equation}\label{ur16.4}
%+\frac{h^2}{6}k_1(x_i)v'''(x_i)+\frac{h^2}{4}k'_1(x_i)v''(x_i)+\frac{h^2}{8}k''_1(x_i)v'(x_i)+\mathcal{O}(h^3).
%\end{equation}
%Similarly we obtain
%\begin{equation}\label{ur16.5}
%\frac{v(x_{i})-v(x_{i-1})}{h}=v'(x_i)-v''(x_i)\frac{h}{2}+v'''(x_i)\frac{h^2}{6}+\mathcal{O}(h^3).
%\end{equation}
%\begin{equation}\label{ur16.6}
%k_1(x_{i-1/2})=k_1(x_i)-k'_1(x_i)\frac{h}{2}+k''_1(x_i)\frac{h^2}{8}+\mathcal{O}(h^3).
%\end{equation}
%From (\ref{ur16.5}) and (\ref{ur16.6}) we get
%$$
%k_1(x_{i-1/2})\frac{v(x_{i})-v(x_{i-1})}{h}=k_1(x_i)v'(x_i)-\frac{h}{2}\left(k_1v'\right)'(x_i)+
%$$
%\begin{equation}\label{ur16.7}
%+\frac{h^2}{6}k_1(x_i)v'''(x_i)+\frac{h^2}{4}k'_1(x_i)v''(x_i)+\frac{h^2}{8}k''_1(x_i)v'(x_i)+\mathcal{O}(h^3).
%\end{equation}
%Subtracting equality  (\ref{ur16.7}) from (\ref{ur16.4}) and
%dividing the result by $h$ we arrive at equality  (\ref{ur16.1}).
%Lemma 5 is proved.

Let  $u(x,t)\in \mathcal{C}_{x,t}^{4,3}$ be a solution of the problem
(\ref{ur01})--(\ref{ur02}).  Let us consider equation  (\ref{ur01})
for $(x,t)=(x_i,t_{j+\sigma})\in\overline Q_T$,\,
$i=1,2,\ldots,N-1$,\, $j=0,1,\ldots,M-1$, $\sigma=1-\alpha/2$:
\begin{equation}\label{ur17}
\partial_{0t_{j+\sigma}}^{\alpha} u=\left.\frac{\partial}{\partial
x}\left(k(x,t)\frac{\partial u}{\partial
x}\right)\right|_{(x_i,t_{j+\sigma})}-q(x_i,t_{j+\sigma})u(x_i,t_{j+\sigma})+f(x_i,t_{j+\sigma}).
\end{equation}

Since
$$
\left.\frac{\partial}{\partial x}\left(k(x,t)\frac{\partial
u}{\partial
x}\right)\right|_{(x_i,t_{j+\sigma})}=k(x_i,t_{j+\sigma})\frac{\partial^2u}{\partial
x^2}(x_i,t_{j+\sigma})+\frac{\partial k}{\partial
x}(x_i,t_{j+\sigma})\frac{\partial u}{\partial x}(x_i,t_{j+\sigma})
$$
$$
=k(x_i,t_{j+\sigma})\left(\sigma\frac{\partial^2u}{\partial
x^2}(x_i,t_{j+1})+(1-\sigma)\frac{\partial^2u}{\partial
x^2}(x_i,t_{j})\right)
$$
$$
+\frac{\partial k}{\partial
x}(x_i,t_{j+\sigma})\left(\sigma\frac{\partial u}{\partial
x}(x_i,t_{j+1})+(1-\sigma)\frac{\partial u}{\partial
x}(x_i,t_{j})\right)+\mathcal{O}(\tau^2)
$$
$$
=\sigma\left(k(x_i,t_{j+\sigma})\frac{\partial^2u}{\partial
x^2}(x_i,t_{j+1})+\frac{\partial k}{\partial
x}(x_i,t_{j+\sigma})\frac{\partial u}{\partial
x}(x_i,t_{j+1})\right)
$$
$$
+(1-\sigma)\left(k(x_i,t_{j+\sigma})\frac{\partial^2u}{\partial
x^2}(x_i,t_{j})+\frac{\partial k}{\partial
x}(x_i,t_{j+\sigma})\frac{\partial u}{\partial
x}(x_i,t_{j})\right)+\mathcal{O}(\tau^2)
$$
$$
=\sigma\left.\frac{\partial}{\partial
x}\left(k(x,t_{j+\sigma})\frac{\partial}{\partial
x}u(x,t_{j+1})\right)\right|_{x=x_i}
$$
$$
+(1-\sigma)\left.\frac{\partial}{\partial
x}\left(k(x,t_{j+\sigma})\frac{\partial}{\partial
x}u(x,t_{j})\right)\right|_{x=x_i}+\mathcal{O}(\tau^2),
$$
$$
q(x_i,t_{j+\sigma})u(x_i,t_{j+\sigma})=q(x_i,t_{j+\sigma})\left(\sigma
u(x_i,t_{j+1})+(1-\sigma)u(x_i,t_{j})\right)+\mathcal{O}(\tau^2),
$$
by virtue of Lemma 5 we have
$$
\left.\mathcal{L}u(x,t)\right|_{(x_i,t_{j+\sigma})}=\sigma\Lambda
u(x_i,t_{j+1})+(1-\sigma)\Lambda
u(x_i,t_{j})+\mathcal{O}(h^2+\tau^2),
$$
where the difference operator $\Lambda$ is defined by formula
(\ref{ur03.5}) with the coefficients
$a_i^{j+1}=k(x_{i-1/2},t_{j+\sigma})$,\,
$d_i^{j+1}=q(x_{i},t_{j+\sigma})$. Let
$\varphi_i^{j+1}=f(x_i,t_{j+\sigma})$, then with regard to Lemma 2
we get the difference scheme with the approximation order
$\mathcal{O}(h^2+\tau^2)$:
\begin{equation}\label{ur18}
\Delta_{0t_{j+\sigma}}^{\alpha}y_i=\Lambda y^{(\sigma)}_i
+\varphi_i^{j+1}, \quad i=1,2,\ldots,N-1,\quad j=0,1,\ldots,M-1,
\end{equation}
\begin{equation}
y(0,t)=0,\quad y(l,t)=0,\quad t\in \overline \omega_{\tau}, \quad
y(x,0)=u_0(x),\quad  x\in \overline \omega_{h},\label{ur19}
\end{equation}

It is interesting to note that for $\alpha\rightarrow1$ we obtain the
Crank--Nicolson difference scheme.

\subsection{Stability and convergence}

\textbf{Theorem 3.} The difference scheme (\ref{ur18})--(\ref{ur19})
is unconditionally stable and its solution satisfies the following a
priori estimate:
\begin{equation}\label{ur20}
 \|y^{j+1}\|_0^2\leq\|y^0\|_0^2+\frac{l^2T^\alpha\Gamma(1-\alpha)}{4c_1}\max\limits_{0\leq j\leq
M}\|\varphi^{j}\|_0^2.
\end{equation}

\textbf{Proof.}   For the difference operator $\Lambda$ using
Green's first difference formula and the embedding theorem
\cite{SamAnd} for the functions vanishing at  $x=0$ and $x=l$, we
get   $(-\Lambda y,y)\geq \frac{4c_1}{l^2}\|y\|_0^2$, that is for
this operator it is possible to take
$\varkappa=\frac{4c_1}{l^2}$.

Since difference scheme (\ref{ur18})--(\ref{ur19}) has the form
(\ref{ur03})--(\ref{ur03.1}), where
$g_s^{j+1}=\frac{c_{j-s}^{(\alpha,\beta)}}{\tau^{\alpha}\Gamma(2-\alpha)}$,
then lemma 4 implies validity of the following inequalities:
$$
g_0^{j+1}=\frac{c_{j}^{(\alpha,\beta)}}{\tau^{\alpha}\Gamma(2-\alpha)}>\frac{1}{2t_{j+\sigma}^\alpha\Gamma(1-\alpha)}>
\frac{1}{2T^\alpha\Gamma(1-\alpha)},
$$
$$
g_j^{j+1}>g_{j-1}^{j+1}>...>g_0^{j+1},
$$
$$
\frac{g_j^{j+1}}{2g_j^{j+1}-g_{j-1}^{j+1}}<\sigma<1.
$$
Therefore, validity of theorem 3 follows from theorem 1. Theorem 3
is proved.

From theorem 2 it follows that if the solution and input data of
problem (\ref{ur01})--(\ref{ur02}) are sufficiently smooth,
the solution of difference scheme  (\ref{ur18})--(\ref{ur19})
converges to the solution of the differential problem with the rate
equal to the order of the approximation error
$\mathcal{O}(h^2+\tau^2)$.

\subsection{Numerical results}

Numerical calculations are performed for a test problem when the
function
$$u(x,t)=\sin(\pi x)\left(t^3+3t^2+1\right)$$
is the exact solution of the problem (\ref{ur01})--(\ref{ur02}) with
the coefficients $k(x,t)=2-\sin(xt)$, $q(x,t)=1-\cos(xt)$ and $l=1$,
$T=1$.

The errors ($z=y-u$) and convergence order (CO) in the norms
$\|\cdot\|_0$ and $\|\cdot\|_{\mathcal{C}(\bar\omega_{h\tau})}$,
where
$\|y\|_{\mathcal{C}(\bar\omega_{h\tau})}=\max\limits_{(x_i,t_j)\in\bar\omega_{h\tau}}|y|$,
are given in Table 2.

\textbf{Table 2} shows that as the number of the spatial
subintervals and time steps is increased keeping $h=\tau$, a
reduction in the maximum error takes place, as expected and the
convergence order of the approximate scheme is
$\mathcal{O}(h^2)=\mathcal{O}(\tau^2)$, where the convergence order
is given by the formula:
CO$=\log_{\frac{h_1}{h_2}}{\frac{\|z_1\|}{\|z_2\|}}$ ($z_{i}$ is the
error corresponding to $h_{i}$).

\textbf{Table 3} shows that if $h=1/1000$, then as the number of
time steps of our approximate scheme is increased, a reduction in
the maximum error takes place, as expected and the convergence order
of time is $\mathcal{O}(\tau^2)$, where the convergence order is
given by the following formula:
CO$=\log_{\frac{\tau_1}{\tau_2}}{\frac{\|z_1\|}{\|z_2\|}}$.

\begin{tabular}{lc}
{\bf Table 2.}\\
$L_2$ - norm and maximum norm error behavior versus grid size reduction \\ when   $\tau=h$.\\
 \hline
$\alpha$ \hspace{12mm} $h$ \hspace{12mm}{$\max\limits_{0\leq n\leq M}\|z^n\|_0$} \hspace{10mm}{CO in $\|\cdot\|_0$} \hspace{8mm}{$\|z\|_{C(\bar \omega_{h \tau})}$} \hspace{8mm}{CO in $||\cdot||_{C(\bar \omega_{h \tau})}$} \\
\hline
0.10      \hspace{3mm}  1/160 \hspace{8mm} $1.0224e-4$ \hspace{11mm}           \hspace{24mm} $1.4518e-4$    \hspace{10mm}         \\
         \hspace{12mm} 1/320 \hspace{8mm} $2.5558e-5$ \hspace{11mm}  2.0001   \hspace{11mm} $3.6294e-5$    \hspace{10mm}   2.0000 \\
         \hspace{12mm} 1/640 \hspace{8mm} $6.3894e-6$ \hspace{11mm}  2.0000   \hspace{11mm} $9.0733e-6$    \hspace{10mm}   2.0000  \\
         \\
0.50      \hspace{3mm}  1/160 \hspace{8mm} $7.8417e-5$ \hspace{11mm}           \hspace{24mm} $1.1153e-4$    \hspace{10mm}         \\
         \hspace{12mm} 1/320 \hspace{8mm} $1.9604e-5$ \hspace{11mm}  2.0000   \hspace{11mm} $2.7882e-5$    \hspace{10mm}   2.0000 \\
         \hspace{12mm} 1/640 \hspace{8mm} $4.9009e-6$ \hspace{11mm}  2.0000   \hspace{11mm} $6.9705e-6$    \hspace{10mm}   2.0000  \\
         \\
0.90      \hspace{3mm}  1/160 \hspace{8mm} $6.6666e-5$ \hspace{11mm}           \hspace{24mm} $9.4949e-5$    \hspace{10mm}         \\
         \hspace{12mm} 1/320 \hspace{8mm} $1.6669e-5$ \hspace{11mm}  1.9998   \hspace{11mm} $2.3740e-5$    \hspace{10mm}   1.9999 \\
         \hspace{12mm} 1/640 \hspace{8mm} $4.1678e-6$ \hspace{11mm}  1.9998   \hspace{11mm} $5.9360e-6$    \hspace{10mm}   1.9998  \\
         \\
0.99     \hspace{3mm}  1/160 \hspace{8mm} $6.5660e-5$ \hspace{11mm}           \hspace{24mm} $9.3532e-5$    \hspace{10mm}         \\
         \hspace{12mm} 1/320 \hspace{8mm} $1.6415e-5$ \hspace{11mm}  2.0000   \hspace{11mm} $2.3384e-5$    \hspace{10mm}   1.9999 \\
         \hspace{12mm} 1/640 \hspace{8mm} $4.1039e-6$ \hspace{11mm}  1.9999   \hspace{11mm} $5.8460e-6$    \hspace{10mm}   2.0000  \\
\hline
\end{tabular}

\vspace{5mm}

\begin{tabular}{lc}
{\bf Table 3.}\\
$L_2$ - norm and maximum norm error behavior versus $\tau$-grid size reduction \\ when $h=1/1000$.\\
 \hline
$\alpha$ \hspace{12mm} $\tau$ \hspace{12mm}{$\max\limits_{0\leq n\leq M}\|z^n\|_0$} \hspace{10mm}{CO in $\|\cdot\|_0$} \hspace{8mm}{$\|z\|_{C(\bar \omega_{h \tau})}$} \hspace{8mm}{CO in $||\cdot||_{C(\bar \omega_{h \tau})}$} \\
\hline
0.10      \hspace{5mm} 1/10 \hspace{8mm} $1.9062e-3$ \hspace{11mm}           \hspace{24mm} $2.6962e-3$    \hspace{10mm}         \\
         \hspace{14mm} 1/20 \hspace{8mm} $4.7789e-4$ \hspace{11mm}  1.9959   \hspace{11mm} $6.7593e-4$    \hspace{10mm}   1.9960 \\
         \hspace{14mm} 1/40 \hspace{8mm} $1.1779e-4$ \hspace{11mm}  2.0205   \hspace{11mm} $1.6659e-4$    \hspace{10mm}   2.0206  \\
         \\
0.50      \hspace{5mm} 1/10 \hspace{8mm} $7.6326e-3$ \hspace{11mm}           \hspace{24mm} $1.0795e-2$    \hspace{10mm}         \\
         \hspace{14mm} 1/20 \hspace{8mm} $1.9130e-3$ \hspace{11mm}  1.9963   \hspace{11mm} $2.7058e-3$    \hspace{10mm}   1.9962 \\
         \hspace{14mm} 1/40 \hspace{8mm} $4.7697e-4$ \hspace{11mm}  2.0039   \hspace{11mm} $6.7461e-4$    \hspace{10mm}   2.0039  \\
         \\
0.90      \hspace{5mm} 1/10 \hspace{8mm} $1.0286e-2$ \hspace{11mm}           \hspace{24mm} $1.4547e-2$    \hspace{10mm}           \\
         \hspace{14mm} 1/20 \hspace{8mm} $2.5706e-3$ \hspace{11mm}  2.0005   \hspace{11mm} $3.6357e-3$    \hspace{10mm}   2.0004  \\
         \hspace{14mm} 1/40 \hspace{8mm} $6.4066e-4$ \hspace{11mm}  2.0045   \hspace{11mm} $9.0608e-4$    \hspace{10mm}   2.0045   \\
         \\
0.99     \hspace{5mm}  1/10 \hspace{8mm} $1.0449e-2$ \hspace{11mm}           \hspace{24mm} $1.4777e-2$    \hspace{10mm}         \\
         \hspace{14mm} 1/20 \hspace{8mm} $2.6102e-3$ \hspace{11mm}  2.0011   \hspace{11mm} $3.6915e-3$    \hspace{10mm}   2.0011 \\
         \hspace{14mm} 1/40 \hspace{8mm} $6.5050e-4$ \hspace{11mm}  2.0045   \hspace{11mm} $9.1998e-4$    \hspace{10mm}   2.0045  \\
\hline
\end{tabular}

\section{A higher order difference scheme for the time
fractional diffusion equation}

%A difference scheme of higher order approximation for the time
%fractional diffusion equation. Stability and convergence. Numerical
%results. }

In this section for problem (\ref{ur01})--(\ref{ur02}), we construct
a difference scheme with the approximation order
$\mathcal{O}(h^4+\tau^2)$ in the case when  $k=k(t)$ and $q=q(t)$.
The stability and convergence of the constructed difference scheme
 in the mesh  $L_2$ - norm with the rate equal to the order of the
  approximation error are proved. The obtained results are
supported by the numerical calculations carried out for a test
example.

\subsection{Derivation  of the difference scheme}

Let us assign a difference scheme to differential problem
(\ref{ur01})--(\ref{ur02}) in the case when  $k=k(t)$ and $q=q(t)$:
\begin{equation}\label{ur21}
\Delta_{0t_{j+\sigma}}^{\alpha}\mathcal{H}_hy_i=a^{j+1}y_{\bar
xx,i}^{(\sigma)}
-d^{j+1}\mathcal{H}_hy_i^{(\sigma)}+\mathcal{H}_h\varphi_i^{j+1}, \,
i=1,\ldots,N-1,\, j=0,1,\ldots,M-1,
\end{equation}
\begin{equation}
y(0,t)=0,\quad y(l,t)=0,\quad t\in \overline \omega_{\tau}, \quad
y(x,0)=u_0(x),\quad  x\in \overline \omega_{h},\label{ur22}
\end{equation}
where $\mathcal{H}_hv_i=v_i+h^2v_{\bar xx,i}/12$, $i=1,\ldots,N-1$,
$a^{j+1}=k(t_{j+\sigma})$, $d^{j+1}=q(t_{j+\sigma})$,
$\varphi_i^{j+1}=f(x_i,t_{j+\sigma})$, $\sigma=1-\alpha/2$.

From  \cite{Sun} and Lemma 2  it follows that if $u\in
\mathcal{C}_{x,t}^{6,3}$, then the difference scheme has the
approximation order  $\mathcal{O}(\tau^2+h^4)$.

\subsection{Stability and convergence}

The difference scheme   (\ref{ur21})--(\ref{ur22}) differs from
(\ref{ur03})--(\ref{ur03.1}) due to the presence of the operator
$\mathcal{H}_h$. However, deriving an a priori estimate for the
solution of difference scheme (\ref{ur21})--(\ref{ur22}) does not
differ significantly from proving Theorem 1.

\textbf{Theorem 4.} The difference scheme (\ref{ur21})--(\ref{ur22})
is unconditionally stable and its solution satisfies the following a
priori estimate:
\begin{equation}\label{ur23}
 \|\mathcal{H}_hy^{j+1}\|_0^2\leq\|\mathcal{H}_hy^0\|_0^2+\frac{l^2T^\alpha\Gamma(1-\alpha)}{c_1}\max\limits_{0\leq j\leq
M}\|\mathcal{H}_h\varphi^{j}\|_0^2,
\end{equation}

\textbf{Proof.} Taking the inner product of the equation
(\ref{ur21}) with
$\mathcal{H}_hy^{(\sigma)}=(\mathcal{H}_hy)^{(\sigma)}$, we have
$$
(\mathcal{H}_hy^{(\sigma)},\Delta_{0t_{j+\sigma}}^{\alpha}\mathcal{H}_hy)-a^{j+1}(\mathcal{H}_hy^{(\sigma)},y_{\bar
xx}^{(\sigma)})
$$
\begin{equation}\label{ur24}
+d^{j+1}(\mathcal{H}_hy^{(\sigma)},\mathcal{H}_hy^{(\sigma)})=(\mathcal{H}_hy^{(\sigma)},\mathcal{H}_h\varphi^{j+1}),
\end{equation}

Let us transform the terms in identity (\ref{ur24}) as
$$
(\mathcal{H}_hy^{(\sigma)},\Delta_{0t_{j+\sigma}}^{\alpha}\mathcal{H}_hy)\geq\frac{1}{2}\Delta_{0t_{j+\sigma}}^{\alpha}\|\mathcal{H}_hy\|_0^2,
$$
$$
-(\mathcal{H}_hy^{(\sigma)},y_{\bar
xx}^{(\sigma)})=-(y^{(\sigma)},y_{\bar
xx}^{(\sigma)})-\frac{h^2}{12}\|y_{\bar
xx}^{(\sigma)}\|_0^2=\|y_{\bar
x}^{(\sigma)}]|_0^2-\frac{1}{12}\sum\limits_{i=1}^{N-1}(y_{\bar
x,i+1}^{(\sigma)}-y_{\bar x,i}^{(\sigma)})^2h
$$
$$
\geq\|y_{\bar x}^{(\sigma)}]|_0^2-\frac{1}{3}\|y_{\bar
x}^{(\sigma)}]|_0^2=\frac{2}{3}\|y_{\bar
x}^{(\sigma)}]|_0^2\geq\frac{8}{3l^2}\|y^{(\sigma)}\|_0^2,\quad
\text{where}\quad \|y]|_0^2=\sum\limits_{i=1}^{N}y_i^2h,
$$
$$
(\mathcal{H}_hy^{(\sigma)},\mathcal{H}_h\varphi^{j+1})\leq\varepsilon\|\mathcal{H}_hy^{(\sigma)}\|_0^2+
\frac{1}{4\varepsilon}\|\mathcal{H}_h\varphi^{j+1}\|_0^2
$$
$$
=\varepsilon\sum\limits_{i=1}^{N-1}\left(\frac{y_{i-1}^{(\sigma)}+10y_{i}^{(\sigma)}+y_{i+1}^{(\sigma)}}{12}\right)^2h+
\frac{1}{4\varepsilon}\|\mathcal{H}_h\varphi^{j+1}\|_0^2\leq\varepsilon\|y^{(\sigma)}\|_0^2+\frac{1}{4\varepsilon}\|\mathcal{H}_h\varphi^{j+1}\|_0^2.
$$
Taking into account the above-performed transformations, from
identity (\ref{ur24}) at $\varepsilon=\frac{8c_1}{3l^2}$ one arrives
at the inequality
$$
\Delta_{0t_{j+\sigma}}^{\alpha}\|\mathcal{H}_hy\|_0^2\leq\frac{l^2}{8c_1}\|\mathcal{H}_h\varphi^{j+1}\|_0^2.
$$
The following process is similar to the proof of theorem 1, and it
is omitted.

The norm $\|\mathcal{H}_hy\|_0$ is equivalent to the norm $\|y\|_0$, which follows from the inequalities
$$
\frac{5}{12}\|y\|_0^2\leq\|\mathcal{H}_hy\|_0^2\leq\|y\|_0^2.
$$

Similarly to theorem 2, we obtain the convergence result.

\textbf{Theorem 5.} Assume that $u(x,t)\in\mathcal{C}_{x,t}^{6,3}$
is the solution of the problem (\ref{ur01})--(\ref{ur02}) in the
case $k=k(t)$, $q=q(t)$, and let $\{y_i^j \,|\, 0\leq i\leq N, \,
1\leq j\leq M\}$ be the solution of the difference scheme
(\ref{ur21})--(\ref{ur22}). Then it holds that
$$
\|u(\cdot,t_j)-y^j\|_0\leq C_R\left(\tau^2+h^4\right),\quad 1\leq
j\leq M,
$$
where $C_R$ is a positive constant independent of $\tau$ and $h$.

\subsection{Numerical results}

In this subsection we present a test example for
a numerical investigation of difference scheme
(\ref{ur21})--(\ref{ur22}).

Consider the following problem:
\begin{equation}\label{ur25}
\partial_{0t}^{\alpha}u(x,t)=k(t)\frac{\partial^2u}{\partial x^2}(x,t)-q(t)u(x,t)+f(x,t),\,\, 0<x<1,\,\, 0<t\leq 1,
\end{equation}
\begin{equation}
u(0,t)=0,\quad u(1,t)=0,\quad 0\leq t\leq 1, \quad u(x,0)=0,\quad
0\leq x\leq 1,\label{ur26}
\end{equation}
where $ k(t)=e^t$, \quad $q(t)=1-\sin{(2t)},$
$$
f(x)=\left[\pi^2t^2e^t+t^2(1-\sin{(2t)})+\frac{2t^{2-\alpha}}{\Gamma(3-\alpha)}\right]\sin(\pi
x),
$$
 whose exact analytical solution reads $u(x,t)=t^2\sin(\pi x).$

\textbf{Table 4} presents the $L_2$ - norm, the maximum norm errors
and the temporal convergence order for $\alpha=0.75, 0.85, 0.95$.
Here we can see that the order of convergence in time is two.

\textbf{Table 5} shows that if  $\tau=1/20000$ is kept fixed, while
$h$ varies, then one obtains the expected fourth-order spatial
accuracy.

\textbf{Table 6} shows that as the number of spatial subintervals
and time steps is increased keeping $h^2=\tau$, a reduction in the
maximum error takes place, as expected and the convergence order of
the approximate of the scheme is $\mathcal{O}(h^4)$.

In \textbf{Table 7} for the case  $N=\lceil\sqrt{M}\rceil$ the
maximum error, the convergence order and CPU time (seconds) are
given. For this case we obtain the expected rate of convergence
$\mathcal{O}(\tau^2)$.

\vspace{3mm}

\begin{tabular}{lc}
{\bf Table 4.}\\
$L_2$ - norm and maximum norm error behavior versus $\tau$-grid size reduction \\ when  $h=1/100$.\\
 \hline
$\alpha$ \hspace{10mm} $\tau$ \hspace{12mm}{$\max\limits_{0\leq n\leq M}\|z^n\|_0$} \hspace{8mm}{CO in $\|\cdot\|_0$} \hspace{8mm}{$\|z\|_{C(\bar \omega_{h \tau})}$} \hspace{8mm}{CO in $||\cdot||_{C(\bar \omega_{h \tau})}$} \\
\hline
0.75     \hspace{3mm}  1/10 \hspace{8mm} $1.6336e-3$ \hspace{11mm}           \hspace{24mm} $2.3103e-3$    \hspace{10mm}         \\
         \hspace{12mm} 1/20 \hspace{8mm} $4.0889e-4$ \hspace{11mm}  1.9983   \hspace{11mm} $5.7826e-4$    \hspace{10mm}   1.9983 \\
         \hspace{12mm} 1/40 \hspace{8mm} $1.0229e-4$ \hspace{11mm}  1.9990   \hspace{11mm} $1.4466e-4$    \hspace{10mm}   1.9990  \\
         \hspace{12mm} 1/80 \hspace{8mm} $2.5581e-5$ \hspace{11mm}  1.9995   \hspace{11mm} $3.6177e-5$    \hspace{10mm}   1.9995  \\
         \\
0.85     \hspace{3mm}  1/10 \hspace{8mm} $1.7130e-3$ \hspace{11mm}           \hspace{24mm} $2.4225e-3$    \hspace{10mm}         \\
         \hspace{12mm} 1/20 \hspace{8mm} $4.2856e-4$ \hspace{11mm}  1.9989   \hspace{11mm} $6.0607e-4$    \hspace{10mm}   1.9989 \\
         \hspace{12mm} 1/40 \hspace{8mm} $1.0718e-4$ \hspace{11mm}  1.9994   \hspace{11mm} $1.5158e-4$    \hspace{10mm}   1.9994  \\
         \hspace{12mm} 1/80 \hspace{8mm} $2.6801e-5$ \hspace{11mm}  1.9997   \hspace{11mm} $3.7902e-5$    \hspace{10mm}   1.9997  \\
         \\
0.95     \hspace{3mm}  1/10 \hspace{8mm} $1.7582e-3$ \hspace{11mm}           \hspace{24mm} $2.4865e-3$    \hspace{10mm}         \\
         \hspace{12mm} 1/20 \hspace{8mm} $4.3967e-4$ \hspace{11mm}  1.9996   \hspace{11mm} $6.2179e-4$    \hspace{10mm}   1.9996 \\
         \hspace{12mm} 1/40 \hspace{8mm} $1.0993e-4$ \hspace{11mm}  1.9998   \hspace{11mm} $1.5547e-4$    \hspace{10mm}   1.9998  \\
         \hspace{12mm} 1/80 \hspace{8mm} $2.7484e-5$ \hspace{11mm}  1.9999   \hspace{11mm} $3.8868e-5$    \hspace{10mm}   1.9999  \\
\hline
\end{tabular}

\vspace{3mm}

\begin{tabular}{lc}
{\bf Table 5.}\\
$L_2$ - norm and maximum norm error behavior versus $h$-grid size reduction \\ when $\tau=1/20000$.\\
 \hline
$\alpha$ \hspace{10mm} $h$ \hspace{14mm}{$\max\limits_{0\leq n\leq M}\|z^n\|_0$} \hspace{8mm}{CO in $\|\cdot\|_0$} \hspace{8mm}{$\|z\|_{C(\bar \omega_{h \tau})}$} \hspace{8mm}{CO in $||\cdot||_{C(\bar \omega_{h \tau})}$} \\
\hline
0.10     \hspace{3mm}  1/4 \hspace{12mm}  $1.1004e-3$ \hspace{11mm}           \hspace{24mm} $1.5562e-3$    \hspace{10mm}         \\
         \hspace{12mm} 1/8 \hspace{12mm}  $6.7512e-5$ \hspace{11mm}  4.0267   \hspace{11mm} $9.5476e-5$    \hspace{10mm}   4.0267 \\
         \hspace{12mm} 1/16 \hspace{10mm} $4.2000e-6$ \hspace{11mm}  4.0067   \hspace{11mm} $5.9397e-6$    \hspace{10mm}   4.0067  \\
         \hspace{12mm} 1/32 \hspace{10mm} $2.6213e-7$ \hspace{11mm}  4.0021   \hspace{11mm} $3.7070e-7$    \hspace{10mm}   4.0021  \\
         \\
0.50     \hspace{3mm}  1/4 \hspace{12mm}  $1.0836e-3$ \hspace{11mm}           \hspace{24mm} $1.5325e-3$    \hspace{10mm}         \\
         \hspace{12mm} 1/8 \hspace{12mm}  $6.6485e-5$ \hspace{11mm}  4.0267   \hspace{11mm} $9.4024e-5$    \hspace{10mm}   4.0267 \\
         \hspace{12mm} 1/16 \hspace{10mm} $4.1360e-6$ \hspace{11mm}  4.0067   \hspace{11mm} $5.8491e-6$    \hspace{10mm}   4.0067  \\
         \hspace{12mm} 1/32 \hspace{10mm} $2.5790e-7$ \hspace{11mm}  4.0034   \hspace{11mm} $3.6472e-7$    \hspace{10mm}   4.0034  \\
         \\
0.90     \hspace{3mm}  1/4 \hspace{12mm}  $1.0654e-3$ \hspace{11mm}           \hspace{24mm} $1.5067e-3$    \hspace{10mm}         \\
         \hspace{12mm} 1/8 \hspace{12mm}  $6.5371e-5$ \hspace{11mm}  4.0266   \hspace{11mm} $9.2449e-5$    \hspace{10mm}   4.0266 \\
         \hspace{12mm} 1/16 \hspace{10mm} $4.0665e-6$ \hspace{11mm}  4.0068   \hspace{11mm} $5.7510e-6$    \hspace{10mm}   4.0068  \\
         \hspace{12mm} 1/32 \hspace{10mm} $2.5346e-7$ \hspace{11mm}  4.0040   \hspace{11mm} $3.5844e-7$    \hspace{10mm}   4.0040  \\
\hline
\end{tabular}

\vspace{3mm}

\begin{tabular}{lc}
{\bf Table 6.}\\
$L_2$ - norm and maximum norm error behavior versus grid size reduction \\ when $h^2=\tau$.\\
 \hline
$\alpha$ \hspace{10mm} $h$ \hspace{14mm}{$\max\limits_{0\leq n\leq M}\|z^n\|_0$} \hspace{8mm}{CO in $\|\cdot\|_0$} \hspace{8mm}{$\|z\|_{C(\bar \omega_{h \tau})}$} \hspace{8mm}{CO in $||\cdot||_{C(\bar \omega_{h \tau})}$} \\
\hline
0.10     \hspace{3mm}  1/10 \hspace{10mm}  $2.4349e-5$ \hspace{11mm}           \hspace{24mm} $3.4434e-5$    \hspace{10mm}         \\
         \hspace{12mm} 1/20 \hspace{10mm}  $1.5166e-6$ \hspace{11mm}  4.0049   \hspace{11mm} $2.1448e-6$    \hspace{10mm}   4.0049 \\
         \hspace{12mm} 1/40 \hspace{10mm}  $9.4708e-8$ \hspace{11mm}  4.0012   \hspace{11mm} $1.3394e-7$    \hspace{10mm}   4.0012  \\
         \hspace{12mm} 1/80 \hspace{10mm}  $5.9180e-9$ \hspace{11mm}  4.0003   \hspace{11mm} $8.3693e-9$    \hspace{10mm}   4.0003  \\
         \\
0.50     \hspace{3mm}  1/10 \hspace{10mm}  $1.4211e-5$ \hspace{11mm}           \hspace{24mm} $2.0097e-5$    \hspace{10mm}         \\
         \hspace{12mm} 1/20 \hspace{10mm}  $8.8285e-7$ \hspace{11mm}  4.0087   \hspace{11mm} $1.2485e-6$    \hspace{10mm}   4.0087 \\
         \hspace{12mm} 1/40 \hspace{10mm}  $5.5094e-8$ \hspace{11mm}  4.0022   \hspace{11mm} $7.7914e-8$    \hspace{10mm}   4.0022  \\
         \hspace{12mm} 1/80 \hspace{10mm}  $3.4420e-9$ \hspace{11mm}  4.0006   \hspace{11mm} $4.8677e-9$    \hspace{10mm}   4.0006  \\
         \\
0.90     \hspace{3mm}  1/10 \hspace{10mm}  $1.5119e-5$ \hspace{11mm}           \hspace{24mm} $2.1381e-5$    \hspace{10mm}         \\
         \hspace{12mm} 1/20 \hspace{10mm}  $9.5080e-7$ \hspace{11mm}  3.9910   \hspace{11mm} $1.3446e-6$    \hspace{10mm}   3.9911 \\
         \hspace{12mm} 1/40 \hspace{10mm}  $5.9571e-8$ \hspace{11mm}  3.9965   \hspace{11mm} $8.4247e-8$    \hspace{10mm}   3.9964  \\
         \hspace{12mm} 1/80 \hspace{10mm}  $3.7274e-9$ \hspace{11mm}  3.9984   \hspace{11mm} $5.2714e-9$    \hspace{10mm}   3.9984  \\
\hline
\end{tabular}

\vspace{3mm}

\begin{tabular}{lc}
{\bf Table 7.}\\
Maximum norm error behavior versus grid size  reduction \\  when $N=\lceil\sqrt{M}\rceil$ and CPU time (seconds).\\
 \hline
$\alpha$ \hspace{12mm} $M$ \hspace{16mm} {$\|z\|_{C(\bar \omega_{h \tau})}$} \hspace{10mm}{CO in $||\cdot||_{C(\bar \omega_{h \tau})}$} \hspace{6mm} CPU(s) \\
\hline
0.70     \hspace{7mm}  10        \hspace{13mm} $2.0986e-3$    \hspace{12mm}                \hspace{33mm}      0.0156       \\
         \hspace{16mm} 30        \hspace{13mm} $2.1085e-4$    \hspace{14mm}   2.0916       \hspace{18mm}      0.0468     \\
         \hspace{16mm} 90        \hspace{13mm} $2.3672e-5$    \hspace{14mm}   1.9905       \hspace{18mm}      0.1404      \\
         \hspace{16mm} 270       \hspace{11mm} $2.6359e-6$    \hspace{14mm}   1.9980       \hspace{18mm}      0.5460     \\
         \hspace{16mm} 810       \hspace{11mm} $2.9428e-7$    \hspace{14mm}   1.9956       \hspace{18mm}      3.0108      \\
         \hspace{16mm} 2430      \hspace{9mm}  $3.2802e-8$    \hspace{14mm}   1.9971       \hspace{16mm}     22.2925      \\
         \\
0.80     \hspace{7mm}  10        \hspace{13mm} $2.1403e-3$    \hspace{12mm}                \hspace{33mm}      0.0156       \\
         \hspace{16mm} 30        \hspace{13mm} $2.2690e-4$    \hspace{14mm}   2.0427       \hspace{18mm}      0.0468     \\
         \hspace{16mm} 90        \hspace{13mm} $2.5342e-5$    \hspace{14mm}   1.9953       \hspace{18mm}      0.1716      \\
         \hspace{16mm} 270       \hspace{11mm} $2.8146e-6$    \hspace{14mm}   2.0004       \hspace{18mm}      0.5616     \\
         \hspace{16mm} 810       \hspace{11mm} $3.1383e-7$    \hspace{14mm}   1.9968       \hspace{18mm}      3.2604      \\
         \hspace{16mm} 2430      \hspace{9mm}  $3.4962e-8$    \hspace{14mm}   1.9976       \hspace{16mm}     23.3065      \\
         \\
0.90     \hspace{7mm}  10        \hspace{13mm} $2.2549e-3$    \hspace{12mm}                \hspace{33mm}      0.0156       \\
         \hspace{16mm} 30        \hspace{13mm} $2.4088e-4$    \hspace{14mm}   2.0358       \hspace{18mm}      0.0468     \\
         \hspace{16mm} 90        \hspace{13mm} $2.6745e-5$    \hspace{14mm}   2.0007       \hspace{18mm}      0.1404      \\
         \hspace{16mm} 270       \hspace{11mm} $2.9607e-6$    \hspace{14mm}   2.0033       \hspace{18mm}      0.5460     \\
         \hspace{16mm} 810       \hspace{11mm} $3.2949e-7$    \hspace{14mm}   1.9986       \hspace{18mm}      3.6670      \\
         \hspace{16mm} 2430      \hspace{9mm}  $3.6670e-8$    \hspace{14mm}   1.9985       \hspace{16mm}     22.7605      \\
\hline
\end{tabular}

\section{Conclusion}

In this paper, the stability and convergence of a family of
difference schemes approximating the time fractional diffusion
equation of a general form is studied. Sufficient conditions for the
unconditional stability of such difference schemes are obtained. For
proving the stability of a wide class of difference schemes
approximating the time fractional diffusion equation, it is simple
enough to check the stability conditions obtained in this paper. A
new difference approximation of the Caputo fractional derivative
with the approximation order $\mathcal{O}(\tau^{3-\alpha})$ is
constructed. The basic properties of this difference operator are
investigated. New difference schemes of the second and fourth
approximation order in space and the second approximation order in
time for the time fractional diffusion equation with variable
coefficients are constructed as well. The stability and convergence
of these schemes in the mesh $L_2$ - norm with the rate equal to the
order of the approximation error are proved. The method can be
easily extended to other time fractional partial differential
equations with other boundary conditions.

 Numerical tests completely
confirming the obtained theoretical results are carried out. In all
the calculations MATLAB is used.

\section{Acknowledgment}

This work was supported by the Russian Foundation for Basic Research
(project 14-01-31246).


\begin{thebibliography}{5}

\bibitem{Nakh:03} A. M. Nakhushev, Fractional Calculus and its
Application, {FIZMATLIT}, {Moscow}, {2003} (in Russian).

\bibitem{Old}
K. B. Oldham, J. Spanier, The Fractional Calculus, Academic Press,
New York, 1974.

\bibitem{Podlub:99}
I. Podlubny, Fractional Differential Equations, Academic Press, San
Diego, 1999.

\bibitem{Hilfer:00}
R. Hilfer (Ed.), Applications of Fractional Calculus in Physics,
World Scientific, Singapore, 2000.

\bibitem{Kilbas:06}
A. A. Kilbas, H. M. Srivastava, J. J. Trujillo, Theory and
Applications of Fractional Differential Equation, Elsevier,
Amsterdam, 2006.

\bibitem{Uchaikin:08}
V. V. Uchaikin, Method of Fractional Derivatives, Artishok,
Ul'janovsk, 2008 (in Russian).

\bibitem{Nigma}
R. R. Nigmatullin, Realization of the generalized transfer equation
in a medium with fractal geometry, Physica Status (B): Basic Res.
133 (1) (1986) 425–430.

\bibitem{Chuk}
K. V. Chukbar, Stochastic transport and fractional derivatives, Zh.
Eksp. Teor. Fiz. 108 (1995), 1875- 1884

\bibitem{Sun}
Y. N. Zhang, Z. Z. Sun, H. L. Liao, Finite difference methods for
the time fractional diffusion equation on non-uniform meshs, J.
Comput. Phys. 265 (2014) 195--210.

\bibitem{Sun2}
Z. Z. Sun, X. N. Wu, A fuly discrete difference scheme for a
diffusion-wave system, Appl. Numer. Math. 56 (2006) 193--209.

\bibitem{Lin}
Y. Lin, C. Xu, Finite difference/spectral approximations for the
time-fractional diffusion equation, J. Comput. Phys. 225 (2007)
1553--1552.

\bibitem{Alikh_arxiv3}
A.A. Alikhanov, {Numerical methods of solutions of boundary value
problems for the multi-term variable-distributed order diffusion
equation}, {arXiv preprint arXiv:1311.2035}, {2013}


\bibitem{ShkhTau:06}
 M. Kh. Shkhanukov-Lafishev, F.I. Taukenova, {Difference methods for solving boundary value problems for fractional differential equations},
    {Comput. Math. Math. Phys.} 46(10) (2006)
    {1785--1795}.





\bibitem{Liu:10}
C. Chen, F. Liu, V. Anh,  I. Turner, {Numerical schemes with high
spatial accuracy for a variable-order anomalous subdiffusion
equations,} {SIAM J. Scien. Comput.} 32(4) (2010) 1740--1760.

\bibitem{Alikh:12}
A.A. Alikhanov, {Boundary value problems for the diffusion equation
of the variable order in differential and difference settings},
{Appl. Math. Comput.} {219} (2012) {3938--3946}.

\bibitem{Delic}
{A. Delic , B.S. Jovanovic,} Numerical approximation of an interface
problem for fractional in time diffusion equation, Appl. Math.
Comput.  229 (2014) 467--479.

\bibitem{Sun3}
R. Du, W. R. Cao, Z. Z. Sun, A compact difference scheme for the
fractional diffusion-wave equation, Appl. Math. Model. 34 (2010)
2998--3007.

\bibitem{Sun4}
G. H. Gao, Z. Z. Sun, A compact difference scheme for the fractional
subdiffusion equations, J. Comput. Phys. 230 (2011) 586--595.

\bibitem{Sun5}
Y. N. Zhang, Z. Z. Sun, H. W. Wu, Error estimates of
Crank-Nicolson-type difference schemes for the subdiffusion
equation, SIAM J. Numer. Anal. 49 (2011) 2302--2322.





\bibitem{Lin2}
Y. Lin, X. Li, C. Xu, Finite difference/spectral approximations for
the fractional cable equation, Math. Comput. 80 (2011) 1369--1396.

\bibitem{Xu}
X. Li, C. Xu, A space-time spectral method for the time fractional
diffusion equation, SIAM J. Numer. Anal. 47 (2009) 2108--2131.

\bibitem{Sun6}
G. H. Gao, Z. Z. Sun, H. W. Zhang, A new fractional numerical
differentiation formula to approximate the Caputo fractional
derivative and its applications, J. Comput. Phys. 259 (2014) 33--50.


\bibitem{Alikh:10}
A.A. Alikhanov, {A priori estimates for solutions of boundary value
problems for fractional-order equations}, {Differ. Equ.} {46(5)}
(2010) {660--666}.




\bibitem{SamAnd}
{\em A. A. Samarskii, V. B. Andreev,} Difference Methods for
Elliptic Equation, Nauka, Moscow, 1976. (in Russian)




\end{thebibliography}
\end{document}